\documentclass[11pt]{amsart}
\usepackage{amsmath}
\usepackage{amssymb}
\usepackage{amsthm}
\usepackage[all]{xy}

\newtheorem{proposition}{Proposition}
\newtheorem{lemma}{Lemma}
\theoremstyle{definition}
\newtheorem{example}{Example}

\CompileMatrices

\begin{document}
\title{A note on sub-bundles of vector bundles}
\author{William Crawley-Boevey}
\author{Bernt Tore Jensen}
\begin{abstract}
It is easy to imagine that a subvariety of
a vector bundle, whose intersection with every fibre is
a vector subspace of constant dimension, must necessarily
be a sub-bundle.
We give two examples to show that this is not true,
and several situations in which the implication does hold.
For example it is true if the base is normal and the field
has characteristic zero.
A convenient test is whether or not the intersections
with the fibres are reduced as schemes.
\end{abstract}

\address{William Crawley-Boevey, \newline
Department of Pure Mathematics \newline
University of Leeds \newline LS2 9JT \newline United Kingdom}

\address{Bernt T. Jensen, \newline
Current address:
Departement de mathematiques \newline
Faculte des Sciences \newline
Universite de Sherbrooke \newline
2500 boul. de l'Universite \newline
Sherbrooke (Quebec) 
J1K 2R1 \newline
Canada}

\maketitle

\section*{Introduction}
Let $X$ be a connected algebraic variety over an algebraically closed 
field~$k$.
Recall that a \emph{vector bundle} of rank $n$ on $X$ is a variety $E$
equipped with a morphism $\pi:E\to X$ and the structure of an 
$n$-dimensional
vector space on each fibre $E_x = \pi^{-1}(x)$, satisfying the local 
triviality
condition that each $x\in X$ has an open neighbourhood $U$ such that 
there
is an isomorphism $\phi:\pi^{-1}(U)\to U\times k^n$ making the diagram
\[
\xymatrix{\pi^{-1}(U) \ar^{\phi}[rr] \ar_{\pi}[dr] & &
U\times k^n \ar^{p}[dl] \\ & U & }
\]
commute, where $p$ is the projection, and with $\phi$ inducing
linear maps on the fibres over $U$.

A \emph{sub-bundle} of $E$ is a closed subvariety $F$ of $E$ 
whose
fibres $F_x = F\cap E_x$ are vector subspaces of the $E_x$ for all 
$x\in X$,
and such that $F$, equipped with the map $\pi|_F:F\to X$ and the 
induced vector
space structures on its fibres, is locally trivial, so a vector bundle.

It is easy to imagine that 
\begin{itemize}
\item[(*)]$F$ is a closed subvariety of $E$ and the $F_x$ are
vector subspaces of the $E_x$ of constant dimension for all $x\in X$,
\end{itemize}
implies that $F$ is a sub-bundle of $E$.
We give two examples showing that this implication is not true in 
general, and several situations in which it does hold.
For example it is true if $X$ is normal and $k$
has characteristic zero. A convenient test is whether or not the intersections
with the fibres are reduced as schemes.

\section{Two general situations}

The first situation is well known.

\begin{proposition}
If $F$ satisfies (*), then $F$ is a sub-bundle if and only if $F$ is 
the kernel of a vector bundle homomorphism $\theta:E\to E'$.
\end{proposition}

If $F$ is a sub-bundle, one can take $E'$ to be the quotient bundle 
$E/F$.
For the converse see for example Proposition 1.7.2 in \cite{lepotier}.
% In all the examples we have found in the literature, there is a natural
% choice for $\theta$, and so the arguments are easily completed.

\begin{proposition}
If $F$ satisfies (*), then $F$ is a sub-bundle if and only if for all 
$x\in X$
the natural scheme structure on $F_x$, as an intersection of varieties 
$F\cap E_x$,
or equivalently as a fibre $(\pi|_F)^{-1}(x)$, is reduced.
\end{proposition}

\begin{proof}
Clearly if $F$ is sub-bundle then local triviality
implies that the $F_x$ are reduced.
For the converse, since the statement is local
on the base, we may assume that $X$ is an affine
variety and $E = X\times k^n$ is trivial.
Let $y=(y_1,\dots,y_n)$ be coordinates for $k^n$
and let $f_1,\dots ,f_s\in k[X][y_1, \dots, y_n]$ be a
generating set for the ideal defining $F$.

For fixed $x$, we identify $E_x$ with $k^n$, and then
$f_1,\dots ,f_s$ define $F_x$ as a closed subscheme of $k^n$.
The tangent space of $F_x$ at $0$ is
\[
T_0 F_x = \{ y\in k^n : J(x)y = 0 \},
\]
where $J(x)$ is the $s\times n$ matrix with
\[
J(x)_{ij} = \frac{\partial f_i(x,y)}{\partial y_j}\biggr|_{y=0}.
\]

We show that $F_x = T_0 F_x$.
If $y\in F_x$, then so is $ty$ for all $t\in k$, so $f_i(x,ty)=0$
for all $i$ and $t$.
Expanding in powers of $t$, and taking the linear terms,
we deduce that $J(x)y=0$, so $F_x \subseteq T_0 F_x$.
Now this is an equality since the scheme structure of $F_x$ is assumed
to be reduced, and the variety structure is that of a vector space,
so $F_x$ is smooth, and hence $\dim T_0 F_x = \dim F_x$.

Now Proposition 1 applies, since $F$ is the kernel
of the vector bundle homomorphism
$\theta:X\times k^n\to X\times k^s$,
$\theta(x,y) = (x,J(x)y)$.
\end{proof}

We observe that in the formulation of (*) we could have used the weaker condition that
$F$ is locally closed in $E$.  

\begin{proposition}\label{closed}
If $F$ is a locally closed subvariety of $E$ and the $F_x$ are
vector subspaces of the $E_x$ of constant dimension for all $x\in X$,
then $F$ is closed in $E$.
\end{proposition}

\begin{proof}
Scalar multiplication on the fibres defines a map
$m:k\times E\to E$, and since the fibres of $F$ are
subspaces, we have $m(k\times F)=F$. It follows that
$m(k\times\overline{F})=\overline{F}$, so that
$\overline{F}\cap E_x$ is closed under scalar multiplication
in the fibre $E_x$, and hence must be connected.
Now $F_x$ is an open subset of $\overline{F}\cap E_x$
since $F$ is locally closed in $E$, and a closed
subset since $F_x$ is a vector subspace of $E_x$.
Thus by connectedness, $F_x = \overline{F}\cap E_x$,
so $F=\overline{F}$.
\end{proof}

\section{Normal base in characteristic zero}

We need the following 
form of Zariski's Main Theorem. For a proof see for example
\S III.9 in \cite{mumford}. 

\begin{lemma}
If $X$ is normal and $f:X'\longrightarrow X$ is a birational
morphism with finite fibres, then $f$ is an isomorphism of $X'$ with an open
subset $U\subseteq X$.
\end{lemma}

\begin{lemma}\label{irred}
If $F$ satisfies (*) and $X$ is irreducible, then so is $F$.
\end{lemma}

\begin{proof}
Let $d$ be the dimension of the $F_x$, and let
$Z_1,\dots,Z_m$ be the irreducible components of $F$.
Applying the theorem on upper semicontinuity of fibre
dimensions to $\pi|_{Z_i}$, see \S I.8 Corollary 3 in \cite{mumford},
one sees that the set
\[
G_i = \{ z\in Z_i \mid \dim_z \left(Z_i \cap \pi^{-1}(\pi(z))\right) 
\ge d \}
\]
is closed in $Z_i$, where $\dim_z Y$ denotes the maximum of the
dimensions of the irreducible components of $Y$ which contain $z$.

Since $m(k\times Z_i)$ contains $Z_i$, and its closure is an
irreducible closed subset of $F$, we must have $m(k\times Z_i)=Z_i$.
Thus the intersection of $Z_i$ with any fibre $E_x$ is closed under
scalar multiplication, and hence any irreducible component of this
intersection contains the zero element. Thus if we define
\[
X_i = \{ x\in X \mid \dim \left( Z_i \cap \pi^{-1}(x) \right) \ge d\},
\]
and $s:X\to E$ denotes the zero section, then
\[
X_i = \{ x\in X \mid \dim_{s(x)} \left( Z_i \cap \pi^{-1}(x)\right) \ge 
d\} = s^{-1}(G_i),
\]
so that $X_i$ is closed in $X$.

Now if $x\in X$, then $F_x = \bigcup_{i=1}^m (Z_i \cap \pi^{-1}(x))$,
so $\dim Z_i \cap \pi^{-1}(x)\le d$. Moreover, since $F_x$ is 
irreducible,
equality must hold for some $i$, and then $F_x = Z_i \cap \pi^{-1}(x)$.
It follows that $X = \bigcup_{i=1}^m X_i$, so since $X$ is irreducible,
$X=X_i$ for some $i$, and then $F=Z_i$, so it is irreducible.
\end{proof}

\begin{proposition}
If $F$ satisfies (*), then $F$ is a sub-bundle if the field $k$ has
characteristic zero and $X$ is normal.
\end{proposition}

\begin{proof}
We have a standing assumption that $X$ is connected, so since $X$ is 
normal
it is irreducible. The statement is local on the base, so by passing to
a neighbourhood of $x_0\in X$,
we may assume that $E = X\times k^n$ is trivial.
For each $x\in X$ identify $F_{x}$ as a vector subspace of $k^n$. Fix a complement 
$C$ of $F_{x_0}$ in $k^n$.

Applying the upper semicontinuity theorem to the restriction of $\pi$
to $F \cap (X\times C)$
one sees that $F_x\cap C=\{0\}$ for all $x$ in a neighbourhood of 
$x_0$.
Then, by shrinking $X$ again, we may assume that $F_x\oplus C = k^n$
for all $x\in X$.

Let $\phi:F\to X\times F_{x_0}$ be the restriction to $F$ of the
map $X\times k^n\to X\times F_{x_0}$ given by the projection
perpendicular to $C$. Clearly $\phi$ is
a bijective morphism of varieties. Now $F$ is irreducible by 
Lemma~\ref{irred},
so $\phi$ is a bijective morphism between
irreducible varieties over a field of characteristic zero,
which shows that $\phi$ is birational by Proposition 7.16 in 
\cite{harris}.
Then $\phi$ is an isomorphism by Zariski's Main Theorem.
This shows that $F$ is locally trivial, and hence a sub-bundle.
\end{proof}

\section{Examples}

\begin{example}
Let $X$ be the cuspidal curve $y^2=x^3$, and let
\[
F = \{ (x,y,z,w) \mid y^2=x^3, yz=xw, w^2=xz^2, x^2z=yw \} \subseteq 
X\times k^2.
\]
We need to know that these equations define $F$ with its reduced scheme
structure, so that they induce the natural scheme structure on the
fibres of $F$. For this it suffices to show that
the ideal corresponding to these equations is the kernel of 
the homomorphism $k[x,y,z,w]\to k[t,z]$ sending $x,y,z,w$ to $t^2,t^3,z,tz$ 
respectively. We leave the details as an exercise to the reader.

Now the fibre of $F$ over a point $(x,y)=(t^2,t^3)$ with 
$t\neq 0$
is the line $w=tz$ in $k^2$, while the fibre over the point $(0,0)$ is
the subscheme of $k^2$ defined by the equation $w^2=0$. As a subset 
this
is the line $w=0$, so every fibre of $F$ is a 1-dimensional subspace of 
$k^2$.
Since the fibre of $F$ over $(0,0)$ is not reduced, $F$ cannot be a
sub-bundle of $X\times k^2$.

%It remains to show that $I=(y^2-x^3, w^2-xz^2, yz-xw, x^2z-yw)$
%is a prime ideal in $k[x,y,z,w]$. For this it suffices to show that
%it is the kernel $K$ of the homomorphism $k[x,y,z,w]\to k[t,z]$
%sending $x,y,z,w$ to $t^2,t^3,z,tz$ respectively.
%Clearly $I\subseteq K$. To show equality, we show that by adding
%suitable elements of $I$, we can reduce an arbitrary element
%of $K$ to zero. Since $w^2-xz^2\in I$, we can assume that the
%element has the form $a(x,y,z)+b(x,y,z)w$. Since $yz-xw\in I$,
%we can subtract any multiple of $x$ from $b(x,y,z)$ provided we
%add the corresponding multiple of $yz$ to $a(x,y,z)$. Similarly,
%since $x^2z-yw\in I$, we can subtract any multiple of $y$
%from $b(x,y,z)$. Thus we may suppose that $b(x,y,z)$ only involves
%powers of $z$, say $b(x,y,z)=f(z)$.
%But now, since the element is in $K$, we have
%$a(t^2,t^3,z) + f(z)tz = 0$. Since the left hand term involves
%no monomials which are linear in $t$, we deduce that $f(z)=0$,
%and so $a(t^2,t^3,z)=0$. Expanding in powers of $z$ and using
%that the kernel of the homomorphism $k[x,y]\to k[t]$ is generated
%by $y^2-x^3$, it follows that $a(x,y,z)\in I$, as required.
\end{example}

\begin{example}
Let $k$ be a field of characteristic $p>0$, and let $X=k$, the
affine line. Let
\[
F = \{ (x,y,z) : y^p = xz^p \} \subseteq k\times k^2
\]
Since $y^p-xz^p$ is an irreducible polynomial, the equation defines
$F$ with its reduced scheme structure. 

The fibre over any point $x=\lambda\in k$ is the line $y = 
\lambda^{1/p}z$,
but the scheme structure is given by the equation $y^p = \lambda z^p$,
or equivalently $(y - \lambda^{1/p}z)^p = 0$, so it is not reduced.
Thus $F$ is not a sub-bundle.
\end{example}

\end{document}